\documentclass[11pt, reqno]{amsart}

\usepackage[colorlinks=true,pagebackref,hyperindex]{hyperref}
\usepackage{mathrsfs} 
\usepackage[all]{xy}
\usepackage{color}
\usepackage{amsfonts}
\usepackage{amscd}

\definecolor{darkgreen}{rgb}{0.0, 0.7, 0.0}

\definecolor{cyan}{cmyk}{1,0,0,0}

\newtheorem{??}[tm]{Question}

\newcommand{\ben}{\begin{enumerate}}
\newcommand{\een}{\end{enumerate}}
\newcommand{\bit}{\begin{itemize}}
\newcommand{\eit}{\end{itemize}}
\newcommand{\beq}{\begin{equation}}
\newcommand{\eeq}{\end{equation}}
\newcommand{\la}{\label}

\newcommand\ci{\cite}

\font\tenmsb=msbm10
\font\sevenmsb=msbm7
\font\fivemsb=msbm5
\newfam\msbfam
\textfont\msbfam=\tenmsb
\scriptfont\msbfam=\sevenmsb
\scriptscriptfont\msbfam=\fivemsb
\def\Bbb#1{{\fam\msbfam #1}}
\font\teneufm=eufm10
\font\seveneufm=eufm7
\font\fiveeufm=eufm5
\newfam\eufmfam
\textfont\eufmfam=\teneufm
\scriptfont\eufmfam=\seveneufm
\scriptscriptfont\eufmfam=\fiveeufm

\newcommand{\im}{ \hbox{\rm Im} }
\newcommand{\ke}{ \hbox{\rm Ker} }
\newcommand{\lorw}{\longrightarrow}

\newcommand\rat{{\Bbb Q}}

\newcommand\comp{{\Bbb C}}
\newcommand\real{{\Bbb R}}
\newcommand\zed{{\Bbb Z}}

\newcommand\pn[1]{{\Bbb P}^{#1}}

\newcommand\s{\sigma}

\newcommand{\m}[1]{\mathcal{#1}}

\newcommand{\td}[1]{ \tau_{ \leq {#1} } }

\newcommand{\Hom}{\mbox{Hom}\,}
\newcommand{\IC}{\m{I}}
\newcommand{\IH}{\mbox{IH}}

\title{The Hodge Theory of Maps\\Lectures 4-5}
\author{
Mark Andrea A.  de Cataldo (Lectures 4-5)}
\thanks{
M.A.  de Cataldo partially supported by N.S.F.}
\date{Sometime in 2010, final revision in 2013}

\begin{document}

\maketitle

\begin{abstract}
These  are the lecture notes
from my two lectures 4 and 5 on the Hodge Theory of Maps delivered at the Hodge Theory Summer School at ICTP Trieste in June 2010. The lectures had a very informal flavor
to them  and, by choice, the notes reflect this fact. They are aimed at beginners, whatever that may mean. There are plenty of exercises and
some references so you can start looking things up on your own.
My little book \cite{dec} contains  some of  the notions discussed here, as well as some amplifications. These notes have appeared as a chapter in the book  {\em Hodge Theory}, Princeton University Press, 2014, edited by E. Cattani, F. El Zein and P. Griffiths. Lecture notes from Lectures 1, 2 and 3 by L. Migliorini  also appear in that book; the present  notes are independent of them, but of course they
are their natural continuation.
\end{abstract}

\tableofcontents

\section{Lecture 4}\label{mdc:lec4}

\subsection{Sheaf cohomology and All That (A Minimalist Approach)}
\label{sc}
\index{cohomology!sheaf cohomology}
\begin{enumerate}

\index{sheaf!injective}\index{injective sheaf}
\item We say that  a sheaf of abelian groups $I$ on a topological space $X$
 is {\em injective} if
\[
\mbox{ the abelian-group-valued functor on sheaves  $\Hom(-,I)$ is exact}.\]
See \ci{ce,gr,ks,ha}.

Of course, the notion of injectivity makes sense in any abelian category,
so we may speak of injective abelian groups, modules over a ring, etc.

\item
{\bf Exercise}

\begin{enumerate}[(a)]

\item
Verify that for every sheaf $F$, the functor $\Hom(-,F)$
is exact on one side (which one?), but, in general,  not on the other.

\item The injectivity of $I$ is equivalent to the following:
for every injection $F\to G$ and every map $F \to I$ there is
a map $G \to I$ making the ``obvious''  diagram (part of the exercise is to
identify this diagram) commutative.

\item
A short exact sequence $0 \to I \to A \to B \to 0$, with $I$ injective, splits.

\item
If $0 \to A \to B \to C \to 0$ is exact and $A$ is injective, then
$B$ is injective if and only if $C$ is.

\item A vector space over a field $k$ is an injective  $k$-module.

\item
By reversing the arrows, you can define the notion of {\em projectivity}
(for sheaves, modules over a ring, etc.). Show that free implies projective.
\index{sheaf!projective}\index{projective sheaf}
\end{enumerate}

\item
\la{z3}
It is a fact that every abelian group can be embedded into an injective
abelian group. Obviously, this is true in the category
of vector spaces!

\item \la{z4}
{\bf Exercise}

 Deduce from the embedding statement above that
every sheaf $F$ can be embedded into an injective sheaf.
(Hint: consider the  direct product sheaf $\Pi_{x \in X} F_x$ on $X$
and work stalk by stalk using \ref{z3}.)

\item
By iteration of the embedding result established in Exercise~\ref{z4},\index{sheaf!injective resolution}
it is easy to show that given every sheaf $F$, there is an
{\em injective resolution of $F$}, i.e.,
 a long exact sequence
\[
0 \lorw F \stackrel{e}\lorw I^0 \stackrel{d^0}\lorw I^1 \stackrel{d^1}\lorw I^2
\stackrel{d^2}\lorw \cdots\]
such that  each $I$ is injective.

\item
The resolution is not unique, but it is so in the homotopy category.
Let us not worry about this; see \ci{ce}
(part of the work to be done by the young (at heart) reader,
is to dig out the relevant statement from the references given here!).
Under suitable assumptions,
usually automatically verified when working with algebraic varieties,
the injective resolution can be chosen to be bounded, i.e., $I^k=0$, for $k \gg 0$;
see \ci{ks}.

\item
Let $f: X \to Y$ be  a continuous map of topological spaces and $F$ be a sheaf on $X$.

The {\em direct image sheaf} $f_* F$  on $Y$  is the sheaf\index{direct image!sheaf}
\[
Y \stackrel{\text{open}}\supseteq U  \longmapsto F (f^{-1}(U)).\]
You should check that the above definition yields a sheaf, not just a presheaf.
\item
A {\em complex of sheaves} $K$ is  a diagram of sheaves and maps of sheaves:\index{complex of sheaves}
\[
 \cdots \lorw K^i \stackrel{d^i}\lorw K^{i+1} \stackrel{d^{i+1}}\lorw \cdots \]
 with $d^2=0$.

 We have the {\em cohomology sheaves}\index{cohomology sheaves}
 \[
 {\mathcal H}^i (K):= \ke{\,d^i}/\,\im{\, d^{i-1}};\]
 recall that everything is first defined as a presheaf and you must take
 the associated sheaf; the only exception is the kernel: the kernel
 presheaf of a map of sheaves is automatically a sheaf (check this).

 A {\em map of complexes} $f: K \to L$ is a compatible system of maps\index{map!of complexes}
 $f^i: K^i \to L^i$. Compatible means that  the ``obvious" diagrams are commutative.

 There are the induced maps of sheaves  ${\mathcal H}^i (f):\index{map!of sheaves}
 {\mathcal H}^i(K) \to {\mathcal H}^i (L)$ for all $i \in \zed$.

A {\em quasi-isomorphism (qis)} $f: K \to L$ is a map inducing isomorphisms
on all cohomology sheaves.\index{quasi-isomorphism}

The {\em translated complex} $K[l]$ has $(K[l])^i:= K^{l+i}$
with the same differentials (up to the  sign $(-1)^l$).\index{translated complex}

Note that $K[1]$ means moving the entries one step to the
{\em left} and taking $-d$.

An exact sequence of complexes is the ``obvious" thing
(make this explicit).\index{exact sequence!of complexes}

Later, I will mention {\em distinguished triangles}:\index{distinguished triangle}
\[ K \lorw  L \lorw M \stackrel{+}\lorw  K[1].\]
You can mentally replace this with a short exact sequence
\[
0 \lorw K \lorw L \lorw M \lorw 0\]
and this turns out to be ok.
\item
The {\em direct image complex} $Rf_* F$ associated with $(F,f)$ is\index{direct image!complex}
``the" complex  of sheaves on $Y$
\[ Rf_* F:= f_* I,\]
where $F \to I$ is an injective resolution as above.\index{injective resolution}

This is well defined up to unique isomorphism in the homotopy\index{homotopy category}
category. This is easy to verify (check it). For the basic definitions
and a proof of this fact
see \ci{ce} (note that there are no sheaves in this reference, the point is the use of the properties of injective objects).

\item
If $C$ is a  {\em bounded below} complex of sheaves on $X$, i.e.,\index{complex of sheaves!bounded below}
 with ${\mathcal H}^i(K) =0$
for all $i \ll 0$ (and we assume this from now on), then $C$ admits a {\em
bounded below injective resolution},
i.e., a qis  $C \to I$, where each entry $I^j$ is injective, and $I$ is  bounded
below.

Again, this is well defined up to unique isomorphism in the homotopy category.

$Rf_*$ is a ``derived functor."
However, this notion and the proof of this fact require plunging into the derived category, which we do not do in these notes. See  \ci{gm}.

\item
We can thus define the {\em derived direct image complex}\index{derived direct image complex}
of a bounded below complex of sheaves $C$ on $X$ by first choosing
a bounded below injective resolution $C \to I$ and then by
setting
\[
Rf_* C := f_* I;\]
this is  a bounded below
 complex of sheaves on $Y$.

\item
Define the {\em (hyper)cohomology groups of  ($X$ with coefficients\index{hypercohomology groups}\index{cohomology!hypercohomology}
in) $C$}  as follows:
\begin{itemize}
\item Take the unique map $c: X \to p$ (a point).

\item Take the complex of global sections $Rc_* C= c_* I= I(X)$.

\item Set  (the right-hand side, being the cohomology of a complex of abelian groups,
 is an abelian group)
\[
H^i (X,C) : = H^i (I(X)).\]
\end{itemize}

\item \la{z13}
{\bf Exercise}  (As mentioned earlier, from now on complexes
are assumed to be bounded below.)

Use the homotopy statements to formulate and prove
that these groups are well defined (typically, this means unique up to unique isomorphism; make this precise).

\item
The {\em direct image sheaves on $Y$ with respect  to $f$}  of\index{direct image!s@sheaves on $Y$ with respect to $f$}
the bounded below complex  $C$ of sheaves  on $X$  are
\[
R^i f_* C : = {\mathcal H}^i (Rf_* C):= {\mathcal H}^i (f_* I),
\qquad i\in \zed.\]
These are well defined (see Exercise~\ref{z13}).

By boundedness,
they are zero for $i \ll 0$ (depending on $C$).

If $C$ is a sheaf,
then they are zero for $i<0$.

\item
{\bf Exercise}

Observe that if $C=F$ is a sheaf, then $R^0f_* F= f_* F$ (as defined earlier).

Prove that the sheaf $R^i f_* C$ is the sheaf associated with the {\em presheaf}
\beq \la{hdi}
U \longmapsto H^i \left(f^{-1} (U), C\right).\eeq
(See \ci{ha}.) This fact is very important in order to build an intuition
for  higher direct images. You should test it against the examples that come to your mind (including all those appearing in these notes).
\index{direct image!higher direct images}

Note that even if $C$ is  a sheaf, then, in general,  (\ref{hdi}) above  defines
 a presheaf. Give many examples
of this fact.

Recall that  while a presheaf and the associated sheaf can be very different,
they have canonically isomorphic stalks! It follows that
 (\ref{hdi}) can be used
to try and determine the stalks of the higher direct image sheaves. Compute these stalks in many examples.

Remark that  for every $y \in Y$ there is a natural map (it is called the base change map)\index{base change map}
\beq\la{bcm}\index{base change map}
(R^i f_* C)_y \lorw H^i (X_y, C_{|X_y})\eeq
between the stalk of the direct image and the cohomology of the fiber
$X_y:= f^{-1}(y)$.

Give examples where this map is not an isomorphism/injective/surjective.

\item
Given a sheaf $G$ on $Y$, the {\em pull-back} $f^* G$ is the sheaf associated\index{sheaf!pull-back}
with the {\em presheaf} (the limit below is the direct limit
over the directed set of open sets $W \subseteq Y$ containing $f(U)$):
\[
U \longmapsto \lim_{W \supseteq f(U)} G(W).\]
This presheaf is not a sheaf even when $f: X \to Y$ is the obvious map
from a set with two elements to a set with one element (both with the discrete
topology) and $G$ is constant.

The pull-back defined above should not be confused with the
pull-back of a quasi-coherent sheaf with respect to a map of algebraic varieties.
(This is discussed very well in  \ci{ha}).

In \ci{go}, you will find a very beautiful discussion of the
\'etale space associated with a presheaf, and hence with a sheaf.
This is all done in the general context of sheaves of sets; it is very
worthwhile
to study sheaves  of sets, i.e., sheaves  without the additional algebraic  structures (sometimes the additional structure
may  hinder some of the basic principles).

The first, important surprise is that every map of sets yields
a sheaf on the target: the sheaf of the local section of the map.

For example, a local homeomorphism,
which can fail to be surjective (by way of contrast, the \'etale
space of a sheaf of abelian groups on a space always surjects onto the space
due to the obvious fact that we always have the zero section!)
yields a sheaf on the target whose \'etale space is canonically isomorphic
with the domain.

Ask yourself: can I view a $2$:$1$ covering space as a sheaf? Yes, see above.
Can I view the same covering as a sheaf of abelian groups? No, unless
the covering is trivial (a sheaf of abelian groups always has the zero section!).

Whereas the definition of direct image $f_*F$ is easy, the \'etale space
of $f_* F$  may bear very little resemblance to the one of $F$.
On the other hand, while the definition of $f^*G$ is a bit more complicated,
the \'etale space $|f^*G|$ of $f^*G$ is canonically isomorphic
with the fiber product over  $Y$ of $X$ with
 the \'etale
space $|G|$ of $G$:
\[
|f^* G|   = |G| \times_Y X.\]

\item\la{adjpr}
It is a fact that  if $I$ on $X$ is injective, then $f_*I$ on $Y$ is injective.

A nice proof of this fact uses the fact that the pull-back functor
$f^*$ on sheaves is the left adjoint to $f_*$, i.e.,  (cf. \ci{gm})
\[
\Hom (f^*F, G) = \Hom(F, f_* G).\]
\index{adjunction property}

\item
{\bf Exercise}

Use the adjunction property in (\ref{adjpr}) to prove that $I$ injective implies
$f_* I$ injective.

Observe that the converse does not hold.

Observe that if $J$ is injective on $Y$, then, in general, the pull-back $f^*J$ is not injective on $X$.

Find classes of maps $f:X \to Y$ for which $J$ injective on $Y$ implies $f^*J$ injective on $X$.

\item{\bf Exercise}

Use that $f_*$ preserves injectives to deduce that
\[
H^i(X,C) = H^i (Y, Rf_* C).\]

\item
It is a fact that on a good space $X$  the cohomology defined
above with coefficients in  the constant sheaf
$\zed_X$ is the same as the one defined by using  singular and
 Cech cohomologies (see \ci{mu,sp}):
\[
H^i(X, \zed_X) = H^i(X,\zed) = \check{H}^i(X, \zed).\]

\item\label{20b.20}
{\bf Exercise} (For a different perspective on what follows,
see Lecture 1.)

\begin{enumerate}[(a)]

\item
Let $j: \real^n - \{0\} \lorw \real^n$ be the open
immersion. Determine the sheaves $R^q j_* \zed$.

 \item\label{20b.b}
 (This is the very first occurrence of the decomposition theorem (DT) in these notes!)
\index{decomposition theorem}
 Let $X=Y=\comp$, $X^*=Y^*=\comp^*$,
 let $f: \comp \to \comp$ be the holomorphic map $z \mapsto z^2$, and let
 $g: \comp^* \to \comp^*$ be the restriction of $f$ to
 $\comp^*:= \comp - \{0\}$.

 Show that $R^i f_* \zed_X =0$ for all $i >0$. Ditto for $g$.

 Show that there is a split short exact sequence of sheaves
 of vector spaces (if you use $\zed$-coefficients,  there is no splitting)
 \[
 0 \lorw \rat_Y \lorw  f_* \rat_X  \lorw Q \lorw 0\]
 and determine the stalks of $Q$.

 Ditto for  $g$ and observe that what you obtain for $g$ is  the restriction to the open set $Y^*$ of what you obtain for $f$ on $Y$.
 (This is a general fact that you may find in the literature
 as ``the base change theorem holds
 for an open immersion.'')

 The  short exact sequence above, when restricted to $Y^*$,\index{sheaf!locally constant}
  is one of locally constant sheaves (recall that
  a locally constant sheaf  of abelian groups---you can guess the
  definition in the case of sheaves of sets with stalk a fixed
  set---with stalk a group $L$ is a sheaf that is locally isomorphic
  to the constant sheaf with stalk $L$)
 and  the restriction  $Q^*$ of $Q$  to $Y^*$
 is the locally constant sheaf with stalk $\rat$ at a point $y \in Y^*$ endowed with
 the automorphism multiplication by $-1$ (explain what this must
 mean).

 The locally constant sheaf $Q^*$  on $Y^*$ (and thus on the unit circle)
 is a good example
 of a nonconstant sheaf with stalk $\rat$, or $\zed$. Another good example is  the sheaf
 of orientations of a nonorientable manifold: the stalk is a set given by two points; this is not a sheaf of groups! If the manifold is orientable,
 then the choice of an orientation turns the sheaf of sets into
  a locally constant sheaf of abelian groups with stalk
  $\zed/  2\zed$.

  \item\label{ex21c.c}
  Show that on a good  connected space $X$,  a locally constant sheaf $L$\index{local system}
  (we often call such an object a {\em local system})
     yields a  representation of the fundamental group $\pi_1 (X,x)$ in the group
     $A(L_x)$ of
     automorphisms of the stalk $L_x$  at a pre-fixed point $x \in X$, and vice versa. (Hint: consider the quotient $(\tilde{X} \times L_x) /\pi_1 (X,x)$
     under a suitable action; here $\tilde{X}$ is a universal covering of $X$.)

     \item
     Use the principle of analytic continuation and the monodromy theorem\index{monodromy theorem}
     (cf. \ci{6iv}) to prove that every local system
     on a simply connected space is constant (trivial representation).

     \item Give an example of a local system that is not semisimple.

     (The relevant definitions are
     simple $:=$ irreducible $:=$ no nontrivial subobject; semisimple $:=$
     direct sum of simples.)

     (Hint: consider, for example,
      the  standard $2 \times 2$ unipotent matrix.)

       The matrix in the hint given above is the one of the Picard--Lefschetz transformation
     associated with the degeneration of a one-parameter family
     of elliptic plane cubic  curves to a rational  cubic curve with a node;
     in other words it is the monodromy of the associated  nontrivial! fiber bundle
     over a punctured disk with fiber $S^1 \times S^1$.

     \item
Given a  fiber bundle,  e.g., a  smooth proper map (see the Ehresmann fibration lemma,\index{Ehresmann fibration lemma}\index{smooth proper map}
e.g., in
\ci{vo}) $f: X
 \to Y$, with fiber $X_y$, prove that  the direct image sheaf is locally constant\index{direct image!sheaf}
 with typical stalk
 \[
( R^i f_* \zed_X)_y =  H^i (X_y, \zed).\]

\item
Show that the Hopf bundle $h: S^3 \to S^2$, with fiber $S^1$, is not
(isomorphic to) a trivial bundle. Though the bundle is not trivial,\index{Hopf!bundle}
the local systems
$R^i h _* \zed_{S^3}$ are trivial on the simply connected $S^2$.

Do the same for $k: S^1 \times S^3 \to S^2$. Verify that
you can turn the above into a proper holomorphic submersion\index{holomorphic!submersion}
of compact complex manifolds $k: S \to \comp \pn{1}$
(see the Hopf surface in \ci{bpv}).

Show that the Deligne theorem (see Lecture 1) on the degeneration for smooth projective maps\index{smooth projective map}
cannot hold for the Hopf map above. Deduce that this is an example
of a map in complex geometry for which the decomposition theorem\index{decomposition theorem}
(DT) does not hold.

  \item
  Show that if  a map $f$ is proper and with finite fibers (e.g., a finite topological covering,  a branched covering,
  the normalization of a complex space, for example of a curve,
  the embedding of a closed subvariety, etc.), then
  $R^if_* F=0$ for every $i >0$ and every sheaf $F$.

  Give  explicit examples of finite maps  and compute $f_* \zed$\index{finite map}
  in those examples.

\end{enumerate}

\item\label{ex.22i.22}
  {\bf Some examples of maps $f: X\to Y$ to play with
  (some have already appeared above)}

\begin{enumerate}[(a)]

\item
$f: (0,1) \to [0,1]$:

$f_* \zed_X  = Rf_* \zed_X=\zed_Y$;

the
base change map (\ref{bcm}) is zero at $0 \in X$.
\index{base change map}

\item
$f: \Delta^*  \to \Delta$ (immersion of punctured unit disk
 into the unit disk
in $\comp$):

$f_* \zed_X = \zed_Y$;

$R^1f_* \zed_X = \zed_o =H^1 (X, \zed_X)$, $o \in \Delta$ the puncture;

there is a nonsplit exact sequence
\[
0 \lorw  \zed_Y \lorw Rf_* \zed_X \lorw \zed_o [-1] \lorw 0.\]

\item
$f: \Delta \to \Delta$, $z \mapsto z^2$:

$R^0 f_* =f_*$; $R^i f_* =0$;
$f_* = Rf_*;$

the natural short exact sequence
\[
0 \lorw R_Y \lorw f_* R_X \lorw Q(R) \lorw 0\]
does not split for $R=\zed$, but it splits  if $2$ is invertible in $R$.

\item
$f: \Delta^* \to \Delta^*$, $z \mapsto z^2$:

$R^0 f_* =f_*$; $R^i f_* =0$;
$f_* = Rf_*;$

the natural short exact sequence
\[
0 \lorw R_Y \lorw f_* R_X \lorw Q(R) \lorw 0\]
does not split for $R=\zed$, but it splits if $2$ is invertible in  $R$;

the stalk $Q(R)_p$ at $p:=1/4 \in \Delta^*$ (the target)
is a  rank-$1$ free $R$-module
generated by the
equivalence class $[(1,-1)]$  in  $R^2/R= (f_* R_X)_p/(R_Y)_p$, modulo
the equivalence relation $(a,b) \sim (a',b')$ if and only if $(a-a'=b-b')$;
here $(a,b)$ is viewed as a constant  $R$-valued function in the
preimage of a  small connected neighborhood of $p$, this preimage being
the disconnected union of two small connected neighborhoods
of $\pm 1/2 \in \Delta^*$ (the domain);

if we circuit once  (e.g., counterclockwise) the origin of the target $\Delta^*$
starting at $1/4$ and returning to it, then the pair
$(1,-1)$ is turned into the pair $(-1,1)$; this is the monodromy representation
on the stalk $Q(R)_p$;

we see that in order to split $R_Y \to  f_* R_X$, or equivalently,
$f_* R_X \to Q(R)$, we need to  be able to divide by $2$.

In this example and in the previous one, the conclusion of the
decomposition theorem (DT) (see Section~\ref{dt}) holds, provided we use
 coefficients in a field of characteristic $0$.\index{decomposition theorem}

 The DT already fails for integer  and for $\zed/2\zed$ coefficients
in these simple examples.

\item
$f: \real \to  Y:= \real /\sim$,  where $Y$ is obtained by identifying $\pm 1 \in \real$
to one  point $o$
(this can be visualized as the real curve $y^2= x^2-x^3$ inside $\real^2$,
with $o$ the origin):

$f_* \zed_X =Rf_* \zed_X$; $(f_* \zed_X)_o \simeq \zed^2$;

there is the natural nonsplit short exact sequence
\[  0\lorw
\zed_Y \lorw f_* \zed_X \lorw \zed_o \lorw 0.\]

Let $j: U:= Y - \{o\} \to Y$ be the open immersion
$j_* \zed = Rj_* \zed$; $(j_* \zed)_o \simeq \zed^4$;

there is the natural nonsplit short exact sequence
\[  0\lorw
\zed_Y \lorw j_* \zed_U \lorw \zed^3_o \lorw 0;\]
note that there is a natural nonsplit short exact sequence
\beq
\la{bndt}
0 \lorw
f_* \zed_X \lorw j_* \zed_U \lorw \zed_o^2 \lorw 0 .
\eeq

\item
$f: \comp \to Y:= \comp /\sim$,  where $Y$ is obtained by identifying $\pm 1 \in \comp$
to a  point $o$ and let  $j: U=Y- \setminus \{o\} \to Y$
(this can be visualized as the complex curve $y^2= x^2-x^3$ inside $\comp^2$,
with $o$ the origin):

this is analogous to the previous example, but it has an entirely different
flavor:
\beq
\la{rfqa}
Rf_* \zed_X = f_* \zed = j_* \zed _U.\eeq
This is another example where the DT holds (in fact here it holds
with $\zed$-coefficients).\index{decomposition theorem}

\item
$f: S^3 \lorw S^3$, the famous Hopf $S^1$-bundle; it is a map of real algebraic\index{Hopf!bundle}
varieties for which the conclusion of Deligne's theorem Section~\ref{dt}(3) does not hold:
we have the trivial local systems\index{Deligne's theorem}
\[
R^0f_* \zed_X = R^1 f_*\zed_X= \zed_Y, \qquad R^if_* \zed_X =0 \ \forall  i \geq 2\]
and a nonsplit (even if we replace $\zed$ with $\rat$) short exact sequence
\[
0 \lorw \zed_Y \lorw Rf_* \zed_X \lorw \zed_Y [-1] \lorw 0\]
(n.b. if it did split, then the first Betti number $0=b_1(S^3)=1!$).

\item
Consider the  action of the group $\zed$ on
$ X':=\comp^2 - \{(0,0)\}$ given by
$(z,w) \mapsto (2z,2w)$.

There are no fixed points and the (punctured
complex lines through the origin, $w=mz$, are preserved).

One shows that $X:= X'/\zed$ is a compact complex surface (a
Hopf surface, see \ci{bpv}) endowed with a proper\index{Hopf!surface}
holomorphic submersion (i.e., with differential
everywhere
of maximal rank) $f: X \to Y= \comp \pn{1}$.

After dividing by the $\zed$-action,
each line $w=mz$ turns into a compact Riemann surface
of genus $1$, which in turn is the fiber $f^{-1} (m)$.
Of course, $m= \infty$ corresponds to the line $z=0$.

If we take the unit $3$-sphere in $\comp^2$, then, $f_{|S^3}:S^3
\to \comp \pn{1} = S^2$
is the Hopf bundle of the previous example.

There is a natural filtration of $Rf_* \zed_X$:
\[
0 = K^{-1} \subseteq K^0 \subseteq K^1 \subseteq K^2= Rf_* \zed_X\]
into subcomplexes
with
\[
K^0/K^{-1} = \zed_Y, \quad K^1/K^0 = \zed_Y^{2} , \qquad K^2/K_1 = \zed_Y.\]
As in the previous example, we cannot have a splitting
\[
Rf_* \zed_Y  \simeq \zed_Y \oplus \zed_Y^2 [-1] \oplus \zed_Y\]
(not even replacing $\zed$ with $\rat$) in view of the fact that
this would imply that $1= b_1 (X) =2$.

This is an example of a proper holomorphic submersion,
where the fibers and the target are projective varieties, but for which
the conclusion of Deligne's theorem  Section~\ref{dt}(3) does not hold.\index{Deligne's theorem}

\item\label{ex.22i.i}
Let $C \subseteq \comp \pn{2}$ be a nonsingular  complex algebraic curve
(it is also  a compact Riemann surface),
let ${\mathcal U}$ be the universal holomorphic line bundle on $\comp \pn{2}$
(the fiber at a point is naturally the complex line parametrized by the point),
let $X$ be the complex surface total space of the line bundle ${\mathcal U}_{|C}$,
let $Y\subseteq  \comp^3$ be the affine cone over $C$; it is a singular surface with an isolated point at the vertex (origin) $o \in Y$.

The blow-up of $Y$ at the vertex  coincides with $X$ (check this).

Let $f: X \to Y$ be the natural map (it contracts the zero section of $X$).

Let $j: U:= Y-\{o\} \to Y$ be the open immersion.

We have the first (for us) example of the DT for a nonfinite map
(for details see \ci{decmigintform}):\index{decomposition theorem}
\beq\la{ret}
Rf_* \rat \;\simeq \;  \td{1} Rj_* \rat_U \; \oplus \;  \rat_o [-2]\eeq
(given a complex $K$, its standard truncated subcomplex
$\td{i}K$ is the complex $L$ with $L^j = K^j$ for every $j <i$,
$L^i := \ke{\, d_K^i}$, $K^j =0$ for every $j >i$; its most important property
is that it has the same cohomology sheaves ${\mathcal H}^j(L)$
 as $K$ for every $j\leq i$ and ${\mathcal H}^j(L) =0$ for every $j >i$).

 The most important aspect of the splitting (\ref{ret}) is that
 the right-hand side does not contain the symbol $f$ denoting the map!
 This is in striking similarity with (\ref{rfqa}), another
 example of DT.
\index{decomposition theorem}

The  relevant direct image sheaves  for $f$ are
\[
R^0 f_* \rat_X =\rat_Y, \quad R^1f_* \rat_X = \rat^{2g}_o, \quad
R^2 f_* \rat_X = \rat_o\]
(the map is proper, the proper base change theorem holds, see  \ci{6iv}, or \ci{ks},\index{proper map}
so that the base change map (\ref{bcm}) is an iso).\index{proper base change theorem}

The  relevant direct image sheaves  for $j$ are
\[
Rj_* \rat_U = \rat_Y, \quad R^1j_*\rat_U = \rat_o^{2g}, \quad
R^2j_* \rat_U = \rat_o^{2g}, \qquad
R^3 j_* \rat_U = \rat_o;\]
this requires a fair amount of work (as a by-product,
you will appreciate the importance of the base change theorem
for proper maps, which you cannot use here!):

$j_* \rat_U= \rat_Y$  is because $U$ is connected;

the computation on the higher $R^ij_* \rat_U$ boils down
to determining
the groups $H^i(U, \rat_U)$ (see (\ref{hdi}));

on the other hand,  $U \to C$ is the $\comp^*$-bundle of the line bundle
${\mathcal U}_{|C}$ and this calculation is carried out in
\ci{bt} (in fact, it is carried out for the associated oriented
$S^1$-bundle) (be warned that \ci{bt} uses the Leray spectral sequence:\index{Leray spectral sequence}
this is a perfect chance to learn about it  without being overwhelmed
by the indices and by being shown very clearly how everything works;
an alternative without spectral sequences is, for example, any
textbook in algebraic topology covering the Wang sequence (i.e.,
the long exact sequence of an oriented $S^1$-bundle; by the way,
it can be recovered using the Leray spectral sequence!).\index{Wang sequence}

Note that if we replace $\rat$ with $\zed$ we lose the splitting
(\ref{ret}) due to torsion phenomena.

Note that there is a  nonsplit short exact sequence
\[
0 \lorw  \rat_Y \lorw \td{1} Rj_* \rat_U \lorw \rat_o^{2g} \lorw 0.\]
A direct proof that this splitting cannot occur is a bit technical
(omitted).

For us it is important to  note that
$\td{1} Rj_* \rat_U$ is the intersection complex\index{intersection complex}
of $I_Y$ of $Y$ (see Section~\ref{ic}) and intersection complexes
$I_Y$ never split  nontrivially
into a direct sum of  complexes.

\item
$f: X= \comp \times C \to Y$, where $C$ is a compact Riemann
surface as in the previous example, where $Y$ is obtained
from $X$ by identifying $\Gamma:= \{0\}\times C$ to a point  $o \in Y$ and leaving the rest
 of $X$ unchanged. Let $U: = Y - \{o\}  = X -  \Gamma$.

 Note  that $\Gamma$ defines the trivial class in $H^2 (X,\zed)$,
 because you can send it to infinity!, i.e., view it as the boundary
 of $\real^{\geq 0} \times C$.

 The actual generator for $H^2(X,\zed) = H^2(C,\zed)$ is given
 by  the class of a complex line $\comp \times c$, $c \in C$.

 You should contrast what is above with the previous example
 given by the total space of a line bundle
 with negative degree. Of course, here $X$ is the total space of the
 trivial line bundle on $C$.

 The map $f$ is not algebraic, not even holomorphic, in fact $Y$
 is not a complex space.

 The  DT cannot hold  for $f$:
\index{decomposition theorem}
 the relevant cohomology sheaves for $Rf_*$ are
 \[
 f_* \rat_X =\rat_Y, \quad, R^1f_* \rat_X = \rat_o^{2g}, \quad
 R^2 f_* \rat = \rat_o;\]
 the relevant cohomology sheaves  for $\td{1} Rj_*\rat_U$
 are
 \[
 j_* \rat_U = \rat_Y, \qquad R^1j_* \rat_U = \rat_o^{2g+1};\]
 it follows that (\ref{ret}), and hence the DT, do not hold in this case.

 For more details and a discussion relating   the first Chern classes of\index{Chern
classes}
 the trivial and of the negative line bundle to the DT
 see \ci{decmigintform}, which also explains (see also \ci{dec}) how to use Borel--Moore homology
 cycles to describe cohomology, as we
 have suggested above.

\end{enumerate}

\end{enumerate}

\subsection{The Intersection Cohomology Complex}\index{intersection cohomology complex}
\label{ic}
We shall limit ourselves  to define and ``calculate" the intersection complex
$I_X$
of a variety of dimension $d$  with one isolated singularity:\index{isolated singularity}
\[
Y = Y_{\mathrm{reg}} \coprod Y_{\mathrm{sing}}, \qquad U:= Y_{\mathrm{reg}}, \qquad Y_{\mathrm{sing}} = \{p\},\]
\[\xymatrix{
U \ar[r]^{j} & Y & \ar[l]_i p.
}
\]
This is done for ease of exposition only. Of course,
the intersection cohomology complex
$I_Y$, and its variants $I_Y(L)$ with twisted coefficients, can be defined for any variety $Y$, regardless of the singularities.
\index{intersection cohomology complex!with twisted coefficients}
\begin{enumerate}
\item
Recall that given a complex $K$ the $a$th {\em truncated complex}\index{truncated complex}
$\td{a} K$ is the subcomplex $C$ with the following entries:
\[
C^b = K^b \; \forall b <a, \qquad C^a = \ke{\, d^a}, \qquad
C^b =0 \; \forall b>a.\]
The single most important property is that
\[
{\mathcal H}^b (\td{a} K) = {\mathcal H}^b (K) \qquad \forall b \leq a,
\qquad \mbox{ zero  otherwise}.\]

\item
Let $Y$ be as above. Define the {\em intersection cohomology complex}
(with coefficients in $\zed$, for example)
as follows:
\[
I_Y:= \td{d-1} Rj_* \zed_U.\]

\item
\la{tomo}
{\bf Toy model}

What follows is related to Section~\ref{sc}, Exercise~\ref{ex.22i.22}\ref{ex.22i.i}.

Let $Y \subseteq \comp^3$ be the affine cone over an elliptic curve
$E \subseteq \comp \pn{2}$.\index{elliptic curve!affine cone over}

$R^0j_* \zed_U  =\zed_Y$ (recall that we always have $R^0f_* =f_*$).

As to the others we observe that $U$ is the $\comp^*$-bundle
of the  hyperplane line bundle $H$ on $E$, i.e., the one induced\index{hyperplane bundle}
by the hyperplane bundle on $\comp \pn{2}$.
By choosing a metric,  we get the unit sphere (here $S^1$) bundle
$U'$ over $E$.  Note that $U'$ and $U$ have the same homotopy type.
The bundle $U' \to E$ is automatically an oriented $S^1$-bundle.
The associated Euler class $e \in H^2 (E, \zed)$
is the first Chern class $c_1 (H)$.\index{Euler class}

\item
{\bf Exercise}

(You will find all you need in \ci{bt}.) Use the spectral sequence\index{Wang sequence}
 for this oriented bundle (here it is just the Wang sequence) to compute
the groups
\[
H^i (U',\zed) = H^i(U, \zed).\]
{\em Answer}: (Caution: the answer below is for $\rat$-coefficients only!; work this situation out in the case of $\zed$-coefficients and keep track of the torsion.)
\[
H^0(U)= H^0(E), \quad H^1(U)= H^1(E), \quad H^2(U)= H^1(E),
\quad H^3(U) = H^2(E).\]
Deduce that, with $\rat$-coefficients (work out the $\zed$ case as well), we have
that $I_Y$ has only two nonzero cohomology sheaves
\[
{\m{H}}^0 (I_Y) = \rat_Y, \qquad
{\m{H}}^1 (I_Y) = H^1(E)_p \;\;\mbox{(skyscraper at $p$)}.\]\index{sheaf!skyscraper}

\item
{\bf Exercise}

Compute $I_Y$ for $Y =\comp^d$, with  $p$ the origin.

{\em Answer}: $I_Y =\rat_Y$ (here $\zed$-coefficients are ok).

\item
 The above result  is general:

 if $Y$ is nonsingular, then $I_Y = \rat_Y$ ($\zed$ ok);

if  $Y$ is the quotient of a nonsingular variety by a finite group action,
then $I_Y = \rat_Y$ ($\zed$ coefficients, KO!).

\item
Let $L$ be a local system on $U$. Define\index{local system}
\[
I_Y(L): = \td{d-1} Rj_* L.\]
Note that (this is a general fact)
\[
{\mathcal H}^0 (I_Y(L)) = j_* L.\]

\item
Useful notation: $j_! L$ is the sheaf on $Y$ which agrees
with $L$ on $U$ and has stalk zero at $p$.

\item
{\bf Exercise}

\begin{enumerate}[(a)]

\item
Let $C$ be a singular curve. Compute $I_C$.\index{algebraic curve!normalization}

{\em Answer}: let $f: \hat{C} \to C$ be the normalization. Then
$I_C= f_* \zed_{\hat{C}}$.

\item

Let things be as in  Section~\ref{sc}, Exercise~\ref{20b.20}\ref{20b.b}.
Let $L = (f_* \zed_X)_{|Y^*}$ and $M := Q_{|Y^*}$.
Compute
\[
I_{Y} (L), \qquad I_Y (M).\]
\item
Let $U$ be as in the toy model Exercise~\ref{tomo}. Determine $\pi_1(U)$.
Classify local systems of low ranks on $U$.
Find some of their $I_Y(L$)'s.

\item
Let $f: C \to D$ be a branched cover of nonsingular
curves. Let $f^o: C^o \to D^o$ be the corresponding topological covering space,
obtained by removing the branching points and their preimages.

Prove that $L:=f^o_* \rat_{C^o}$ is semisimple ($\zed$-coefficients
is KO!, even for the identity: $\zed$ is not a simple $\zed$-module!).

Determine $I_D(L)$ and describe its stalks. (Try the case when $C$ is replaced by a surface, threefold, etc.)

\end{enumerate}

\end{enumerate}

\subsection{Verdier Duality}\index{Verdier duality}
\label{vd}
For ease of exposition, we work with rational coefficients.

\begin{enumerate}
\item
Let $M^m$ be an oriented manifold. We have Poincar\'e
duality:\index{Poincar\'e duality}
\beq\la{pdeq}
H^i(M,\rat) \simeq H^{m-i}_c (M,\rat)^*.\eeq

\item
{\bf Exercise}

Find compact and noncompact examples
of the failure of Poincar\'e duality for singular complex varieties.

(The easiest way to do this is to  find nonmatching Betti numbers.)

\item
Verdier duality (which we do not define here; see \ci{gm})
is the culmination of a construction that  achieves the following
generalization  of Poincar\'e duality to the case of
complexes of sheaves  on  locally compact spaces.

Given a complex of sheaves $K$ on $Y,$ its  {\em Verdier dual}
$K^*$ is a canonically defined  complex on $Y$  such that
 for every open $U \subseteq Y$,
\beq\la{vdeq}
H^i (U, K^*) = H^{-i}_c (U, K)^*.\eeq
Note that $H_c^i(Y,K)$ is defined the same way as $H^i(Y,K)$, except
that we take global sections with compact supports.

The formation of $K^*$ is contravariantly functorial in $K$:
\[
K \lorw L , \qquad K^* \longleftarrow  L^*,\]
and satisfies
\[
K^{**}=K, \qquad (K[l])^*= K^* [-l].\]

\item
{\bf Exercise}

Recall the definition of the translation\index{translation functor}
functor $[m]$ on complexes (see Section~\ref{sc}) and those of $H^i$ and $H^i_c$
and show  directly that
\[
H^i(Y, K[l]) = H^{i+l}(Y,K), \qquad H^i_c(Y, K[l]) = H^{i+l}_c(Y, K).\]

\item
It is a fact that, for the oriented manifold $M^m$,
the chosen orientation determines an isomorphism
\[
\rat_Y^*= \rat_Y [m]\]
so that we get Poincar\'e duality.
Verify this!; that is, verify that (\ref{vdeq}) $\Longrightarrow$ (\ref{pdeq})

(do not take it for granted, you will see what duality means over a point!).

If $M$ is not oriented, then you get something else. See \ci{bt} (look for ``densities"), see \ci{ks} (look for ``sheaf of orientations"), see \ci{6iv} (look for  ``Borel--Moore
chains"),
and the resulting complex of sheaves (see also \ci{ih}).

\item
One of the most important properties of $I_Y$ is its self-duality, which we express
as follows (the translation by $d$  is for  notational convenience):
first set
\[
\IC_Y: = I_Y[d]\]
(we have  translated the complex $I_Y$, which had  nonzero cohomology sheaves
only in degrees $[0,d-1]$, to the \emph{left} by $d$ units, so that
the corresponding interval is now $[-d,-1]$); then we have that\index{cohomology sheaves}
\[
\IC_Y^*=\IC_Y.\]

\item
{\bf Exercise}

 Use the toy model to verify that the equality
holds (in that case) at the level of cohomology sheaves
by verifying that (here $V$ is a  ``typical" neighborhood of $p$)
\[
{\mathcal H}^i(\IC_Y)_p = H^i(V, \IC_V) =H^{-i}_c (V, \IC_V)^*.\]
(To do this, you will need to compute $H^i_c(U)$ as you did
$H^i(U)$; be careful though about using homotopy
types and $H_c$!)
You will find the following distinguished triangle useful\index{distinguished triangle}---recall
we can view them as short exact sequences, and as such, yielding a long exact sequence of cohomology
groups, with or without supports:
\[
{\mathcal H}^0 (I_Y)  \lorw I_Y \lorw {\mathcal H}^1(I_Y) [-1] \stackrel{+}\lorw;\]
you will also find useful the long exact sequence
\[
\cdots \lorw
H_c^a(U) \lorw H_c^a (Y) \lorw H^a_c (p) \lorw H^{a+1}_c(U) \lorw \cdots.\]

\item
Define the {\em intersection cohomology groups of $Y$} as
\[
\IH^i(Y) = H^i(Y, I_Y), \qquad \IH^i_c(Y) = H^i_c(Y, I_Y).\]
The original definition is more geometric and involves chains
and boundaries, like in the early days of  homology; see \ci{ih}.

\item
Since $\IC_Y^*=\IC_Y$, we get
\[H^i(Y, \IC_Y)= H^{-i}_c(Y, \IC_Y)^*.\]
Using $\IC_Y= I_Y [d]$, Verdier duality implies that
\[
H^i(Y, I_Y) = H^{2n-i}_c(Y, I_Y)^*,
\]
and we immediately deduce  {\em Poincar\'e duality for intersection cohomology
 groups}
on an arbitrarily singular complex algebraic variety (or complex space):
\[
\IH^i(Y, I_Y) = \IH^{2d-i}_c(Y, I_Y)^*.\]

\item
Variant for twisted coefficients.

If $Y^o \subseteq Y_{\mathrm{reg}} \subseteq Y$, $L$ is a local system
on  a nonempty open set $Y^o$ and  $L^*$  is the dual local  system,
then we have $I_Y(L)$, its translated $\IC_Y(L)$, and
we have a canonical isomorphism
\[
\IC_Y(L)^* = \IC_Y(L^*).\]
There is the corresponding duality statement
for the groups $\IH^i(Y, I_Y(L))$, etc.:\[
\IH^i(Y, I_Y(L^*)) = \IH^{2d-i}_c(Y, I_Y(L))^*.\]

\item
{\bf Exercise}

Define the dual local system $L^*$ of a local system $L$ as the
sheaf of germs of sheaf maps $L \to \rat_Y$.\index{local system!dual}

\begin{enumerate}
\item
Show that it is a local system and that
there is a pairing (map of sheaves)\index{map!of sheaves}
\[
L \otimes_{\rat_Y} L^* \lorw \rat_Y\]
inducing identifications
\[
(L_y)^* = (L^*)_y.
\]
(Recall that the tensor product is defined by taking the sheaf associated with
the presheaf tensor product (because of local constancy
of all the players, in this case the presheaf is a sheaf):
$U \mapsto L(U) \otimes_{\rat_U(U)} L^*(U)$).

\item
If $L$  is given by the representation $r: \pi_1 (Y,y) \to A(L_y)$
(see Section~\ref{sc}, Exercise~\ref{20b.20}\ref{ex21c.c}),
find an  expression for a representation  associated with $L^*$.
(Hint: inverse--transpose.)

\end{enumerate}

\item Verdier duality and $Rf_*$ for a proper map.\index{direct@$Rf_*$ for a proper map}\index{proper map}
\label{vdrf}

It is a fact that if $f$ is proper, then
\[
(Rf_* C)^* = Rf_* (C^*).\]
We apply this to $\IC_Y (L)^*=\IC_Y(L^*)$ and get
\[
(Rf_* \IC_Y(L))^*  = Rf_* \IC_Y(L^*).\]
In particular, $Rf_* \IC_Y$  is self-dual.

\end{enumerate}

\section{Lecture 5}
\subsection{The Decomposition Theorem (DT)}
\label{dt}
\index{decomposition theorem}
\begin{enumerate}
\item
Let $f: X \to Y$ be a {\em proper} map of algebraic varieties and $L$\index{local system!semisimple}\index{semisimple local system}
be a {\em semisimple}  ($=$ direct sum of simples; simple $=$ no nontrivial
subobject) local system with $\rat$-coefficients\index{proper map}
(most of what follows fails with coefficients
not in a field of characteristic $0$)

 on a Zariski  dense open
set $X^o \subseteq X_{\mathrm{reg}} \subseteq X$.

Examples include
\begin{itemize}
\item $X$ is nonsingular, $L= \rat_X$, then $I_X (L)= I_X =\rat_X$;

\item $X$ is singular, $L= \rat_{X_{\mathrm{reg}}}$, then $I_X(L) = I_X$.
\end{itemize}

\item
{\bf Decomposition theorem}

The following statement is the deepest known fact
concerning the homology of algebraic varieties.

{\em
There is a splitting in the derived category of sheaves on $Y$:\index{derived category!splitting}
\beq \la{wdt}
Rf_*  \, I_X(L) \simeq   \bigoplus_{b \in B} I_{\overline{Z_b}} (L_b) [l_b],\eeq
where
\begin{itemize}
\item B is a finite set of indices;

\item $Z_b\subseteq Y$ is a collection of locally closed nonsingular subvarieties;

\item $L_b$ is a semisimple local system on $Z_b$; and

\item $l_b \in \zed$.
\end{itemize}
}
\bigskip

What does it mean to have a splitting in the derived category?

Well, I did not define what a derived category is (and I will not).
Still,   we can deduce immediately from (\ref{wdt})  that
the intersection cohomology  groups of the domain split into  a direct sum
of intersection cohomology groups on the target.

\item
The case where we take $I_X= I_X(L)$ is already important.

Even if $X$ and $Y$ are smooth, we must deal with $I_Z$'s on $Y$, i.e.,
we cannot have a direct sum of shifted sheaves for example.

Deligne's theorem (1968), including the semisimplicity statement
(1972) for proper smooth maps of smooth varieties (see Lectures 1 and 2) is\index{smooth proper map}
a special case and it reads as follows:
\[
R f_* \, \rat_X \simeq \bigoplus_{i\geq 0} R^if_* \rat_X [-i],
\qquad I_Y(R^if_* \rat_X) = R^if_* \rat_X.\]

\index{Deligne's theorem}

\item\label{symm}

{\bf Exercise}

By using the self-duality of
$\IC_Y$,  the rule $(K[l])^*=K^*[-l]$, the DT above, and the fact that $\IC_T = I_T [\dim{T}]$, show that (\ref{wdt})
 can be rewritten  in the following  more
symmetric form, where $r$ is a uniquely determined nonnegative integer:
\[
Rf_* \IC_X \simeq
\bigoplus_{i =- r}^r P^i[-i],\]
where  each $P^i$ is a direct sum of some of the $\IC_{\overline{Z_b}}$
appearing above, {\em without translations $[-]$}!, and
\[
(P^i)^* =P^{-i} \qquad \forall i \in \zed.\]
Try this first in the case of smooth proper maps, where
$Rf_*\rat_X= \oplus R^if_*\rat_X [-i]$. This may help  to get used to the change of indexing scheme
as you go from $I_Y$ to $\IC_Y= I_Y[d]$.\index{smooth proper map}

\item\label{stc.5}
{\bf Exercise}

\begin{enumerate}[(a)]
\item
Go back to all the examples we met earlier and determine, in the cases where
the DT is applicable,
the summands appearing on the left of (\ref{wdt}).
\index{decomposition theorem}

\item
\label{stc.b}
(See \ci{decmigintform}.)
Let $f: X \to C$ be a proper algebraic  map with connected fibers, $X$ a nonsingular\index{proper algebraic map}
algebraic
surface, $C$ a nonsingular algebraic curve.

Let $C^o$ be the set of regular values, $\Sigma:= C\setminus C^o$ (it is  a fact that it is finite). Let $f^o: X^o \to C^o$ and $j: C^o \to C$ be the obvious maps.

Deligne's theorem applies to $f^o$ and is a statement on $C^o$;
show that it takes the following form:
\[
Rf^o_* \rat_{X^o} \simeq   \rat_{C^o} \oplus R^if^o_* \rat_{X^o} [-1] \oplus
\rat_{C^o} [-2].\]
Show that the DT on $C$ must take the form (let $R^1:= R^1 f^o_* \rat_{X^o}$)
\[
Rf_* \rat_X \simeq  \rat_C \oplus j_* R^1 [-1] \oplus \rat_C [-2] \oplus V_{\Sigma}[-2],\]
where $V_{\Sigma}$ is the skyscraper sheaf on the finite set\index{sheaf!skyscraper}
$\Sigma$ with stalk at each $\s \in \Sigma$ a vector space $V_\s$
of rank equal to the number of irreducible components of $f^{-1}(\s)$
{\em minus $1$}.

Find a more canonical description of $V_{\s}$ as a quotient of $H^2(f^{-1}(\s))$.

Note that this splitting contains quite a lot of information. Extract it:

\begin{itemize}
\item
  The only feature of $f^{-1} (\s)$ that contributes to $H^*(X)$ is its number of irreducible components; if this is $1$, there is no contribution,
  no matter how singular (including multiplicities) the fiber is.

\item
Let $c \in C$, let   $\Delta$ be a small  disk around $c$,
let $\eta \in \Delta^*$ be a regular value.

We have the bundle $f^* : X_{\Delta^*} \to \Delta^*$ with typical fiber
$X_\eta:= f^{-1} (\eta)$.

 We have the (local) monodromy for this bundle; i.e., $R^i$ is a local system;\index{monodromy}\index{local system}
 i.e.,
  $\pi_1 (\Delta^*) =\zed$ acts on $H^i(X_\eta)$.

  Denote by ${R^1}^{\pi_1} \subseteq R^1_\eta$ the invariants of this (local) action.

  Show the following  general fact: for local systems $L$ on a
  good  connected space $Z$ and for a point $z \in Z$ we have that
  the invariants of the local system $L_z^{\pi_1(Z,z)}= H^0 (Z,L)$.

Let $X_c:= f^{-1}(c)$ be the central fiber;
there are the natural
restriction maps
\[
 H^1(X_{\eta} ) \supseteq H^1 (X_\eta)^{\pi_1} \stackrel{r}\longleftarrow H^1 (f^{-1} (\Delta) )
 \stackrel{\simeq}\lorw H^1(X_{c}).\]

 Use the DT above to deduce that $r$ is {\em surjective}---this
 is the celebrated {\em local invariant cycle theorem}: all local invariant classes
 come from $X_{\Delta}$; it comes {\em for free} from the DT.\index{local invariant cycle theorem}\index{decomposition theorem}

 \end{itemize}

 Finally observe, that in this case, we indeed have
 $Rf_* \rat_X \simeq \oplus R^i f_* \rat [-i]$ (but you should view this as
 a coincidence due to the low dimension).

 \item
 Write down the DT for a projective bundle over a smooth variety.\index{decomposition theorem}

 \item
 Ditto for the blowing up of a nonsingular subvariety of a nonsingular variety.

 \item\label{3f.3}
 Let $Y$ be a threefold with an isolated singularity at $p \in Y$.
 Let $f: X \to Y$ be a resolution of the singularities of $Y$:
 $X$ is nonsingular, $f$ is proper and it is
 an isomorphism over $Y - \{ p\}$.

 \begin{enumerate}
 \item
 \label{3f.f} Assume $\dim{f^{-1} (p)} =2$; show, using
 the symmetries expressed by Exercise~\ref{symm}, that the DT
 takes the form
 \[
 Rf_* \rat_X= \IC_Y \oplus V_p [-2] \oplus W_p[-4],\]
 where $V_p\simeq W_p^*$ are skyscraper sheaves with dual stalks.

Hint: use
 $H_4 (X_p) \neq 0$ (why is this true?) to infer, using that ${\mathcal H}^4
 (I_X)=0$, that one must have a summand contributing to $R^4 f_*\rat$.

 Deduce that the irreducible components of top dimension $2$ of $X_p$
 yield linearly independent cohomology classes in $H^2(X)$.

 \item
 Assume $\dim{ f^{-1} (p)} \leq 1$. Show that
 we must have
 \[
 Rf_* \rat_X = I_Y.\]
 Note that this is remarkable and highlights a general principle:
 the proper algebraic maps are restricted by the fact that the topology of $Y$,\index{proper algebraic map}
  impersonated by $I_Y$, restricts the topology of $X$.

As we have seen in our examples to play with at the end of
Section~\ref{sc}, there are no such general restriction in other geometries,
 e.g., proper $C^\infty$ maps,\index{proper smooth map}
proper  real algebraic maps,   proper holomorphic maps.

 \end{enumerate}

\end{enumerate}

\end{enumerate}

\subsection{(Relative) Hard Lefschetz for Intersection Cohomology}
\label{rhl}

\begin{enumerate}
\index{hard Lefschetz theorem!relative}\index{hard Lefschetz theorem!for intersection cohomology groups}
\item
Let $f: X \to Y$ be a projective smooth map of nonsingular varieties\index{projective smooth map}
and $\ell \in H^2 (X,\rat)$ be the first Chern class of a line bundle
on $X$ which is ample (Hermitian positive) on every $X_y$.

We have the iterated  cup product map (how do you make this precise?)
\[
{\ell}^i: R^j f_* \rat_X \lorw R^{j+2i} f_* \rat_X.\]
For every fiber $X_y:= f^{-1} (y)$,
we have the {\em hard Lefschetz theorem} (\ci{6gh}) for  the  iterated cup product action of
$\ell_y \in H^*(X_y, \rat)$; let $d = \dim{X_y}$.

The hard Lefschetz theorem on the fibers of the smooth proper map\index{smooth proper map}
$f$  implies at once that we have the  isomorphisms of sheaves
\beq \la{rhls}
\ell^i: R^{d-i} f_* \rat_X \stackrel{\simeq}\lorw R^{d+i}f_* \rat_X\eeq
and we view this fact  as the {\em relative hard Lefschetz
theorem for smooth proper maps}.
\index{relative hard Lefschetz theorem}
\item
In an earlier exercise, you were asked to find examples of the failure of
Poincar\'e duality. It was suggested you find examples of
(necessarily singular)
complex projective varieties of complex dimension $d$ for which one does not have
the  symmetry predicted by Poincar\'e duality: $b_{d-i} = b_{d+i}$, for every  $i \in \zed$.
\index{Poincar\'e duality}
Since the conclusion of the hard Lefschetz theorem yields the same symmetry
for the Betti numbers, we see that for these same examples, the conclusion
of the hard Lefschetz theorem does not hold.\index{Betti numbers}

If the hard Lefschetz theorem does not hold for singular projective varieties,
the sheaf-theoretic counterpart (\ref{rhls}) cannot hold (why?) for an arbitrary
proper map, even if the domain and target are nonsingular and the map
is surjective (this is due to the singularities of the fibers).\index{hard Lefschetz theorem}\index{proper map}

 In short, the relative hard Lefschetz does not hold if formulated in terms
 of  an isomorphism between  direct image {\em sheaves}.

 \item
 Recall the symmetric form of the DT (see Section~\ref{dt}, Exercise~\ref{symm}):
 \[Rf_* \IC_X \simeq \bigoplus_{i=-r}^r P^i[-i].\]
\index{decomposition theorem}

 It is a formality
 to show that given a map $f: X \to Y$ and a cohomology class $\ell \in H^2 (X, \rat)$
we get iterated cup product maps\index{cup product}
\[
\ell^i : P^j \to P^{j+2i}.\]

 The {\em relative hard Lefschetz theorem} (RHL) is the statement that
 if $f$ is proper and
 if $\ell$ is the first Chern class of an ample line bundle on $X$,\index{Chern
classes}
 or at least ample on every fiber of $f$, then we have that the iterated cup product maps
 \beq\la{rrhhll} \ell^i: P^{-i} \stackrel{\simeq}\lorw  P^i\eeq
 are isomorphisms for every $i\geq 0$.

 In other words, the conclusion of the RHL
 (\ref{rhls})
 for smooth proper maps, expressed as an isomorphism of direct image sheaves,
 remains valid for arbitrary proper maps provided\index{smooth proper map}

 \begin{itemize}
 \item we push forward $\IC_X$, i.e., we form $Rf_* \IC_X$, vs.~$Rf_* \rat_X$
 for which nothing so clean holds in general; and

 \item we consider the complexes $P^i$, instead of the direct image sheaves.

 \end{itemize}

In the interest of perspective, let me add that the $P^i$
 are the so-called perverse direct image complexes of $\IC_X$
 with respect to $f$ and are special
 perverse sheaves on $Y$. The circle of ideas is now closed:

{\em  \centerline{ RHL is a statement about the perverse direct image
 complexes of $Rf_* \IC_X$!}}
\index{perverse direct image complex}

 Note that Verdier duality shows that $P^{-i} = (P^i)^*.$
 Verdier duality holds in general, outside of the realm of
 algebraic geometry
 and holds, for example for the Hopf surface map $h: S \to \comp \pn{1}$.
\index{Verdier duality}
 In the context of complex geometry,
 the RHL, $\ell^i: P^{-i} \simeq P^i$, is a considerably deeper statement than Poincar\'e duality.

 \item
 {\bf Exercise}
 \begin{enumerate}[(a)]
 \item
 Make the statement of the RHL explicit in the example
 of a map from a surface to a curve (see Section~\ref{dt}, Exercise~\ref{stc.5}\ref{stc.b}).

 \item
 Ditto for  Section~\ref{dt}, Exercise~\ref{stc.5}\ref{3f.3}\ref{3f.f}.
 (Hint: in this case you get $\ell: V_p \simeq W_p$.)

 Interpret geometrically, i.e., in terms of intersection theory,
 the isomorphism $i: V_p \simeq W_p^*$ (PD) and $l: V_p \simeq W_p$ (RHL).

 (\emph{Answer}: (see \ci{decmigintform}) let $D_k$ be the fundamental classes of the exceptional
 divisors (which are the surfaces in $X$ contracted to $p$); interpret
 $W_p$ as (equivalence classes of) topological  $2$-cycles $w$; then $i$ sends
 $D_k$ to the linear map sending $w$ to $D_k \cdot w \in H^6(X,\rat) \simeq \rat$;
 the map  should be viewed as the operation of intersecting
 with a  hyperplane section $H$  on $X$ and it sends
 $D_k$  to the $2$-cycle $D_k \cap H$. Now you can word out the conclusions
 of PD and RHL and appreciate them.)

 \end{enumerate}

 \item
 The hard Lefschetz theorem on  the intersection cohomology
 groups $\IH(Y,\rat)$
 of a projective variety $X$ of dimension $d$.
\index{hard Lefschetz theorem!for intersection cohomology groups}
 Let us apply RHL to the  proper map
 $X \to \mathrm{point}$:

 let $\ell$ be the first Chern class of an ample line bundle
on  $X$ of dimension $d$, then
 \[\ell^i: \IH^{d-i}(X, \rat)  \stackrel{\simeq}\lorw \IH^{d+i}(X,\rat).\]

 \item
 Hodge--Lefschetz package for intersection cohomology.\index{Hodge Lefschetz package@Hodge--Lefschetz package!for intersection cohomology}

 Let $X$ be a projective variety. Then the statements
 (see \ci{6gh} for these statements)
 of the two (hard and hyperplane section) Lefschetz theorems,
 of the primitive Lefschetz decomposition,  of\index{Lefschetz decomposition!for intersection cohomology}
 the Hodge decomposition and of  the Hodge--Riemann bilinear relations\index{Hodge Riemann bilinear relations@Hodge--Riemann bilinear relations!for intersection cohomology}
 hold for the rational intersection cohomology group of  $\IH(X, \rat)$.

 \item
 {\bf Exercise} (Compare  what follows with the first part of Lecture 3.)

 Let $f: X \to Y$ be a resolution of the singularities of a projective surface with isolated singularities  (for simplicity only; after you solve this, you may want to tackle the case when the singularities are not isolated).

Show that the DT takes the form\index{decomposition theorem}
\[
Rf_* \rat_X[2] = \IC_Y \oplus V_{\Sigma},\]
where
$\Sigma$ is the set of singularities of $Y$ and $V_{\Sigma}$ is the skyscraper
with fiber $V_\s = H^2(X_\s)$ (here $X_\s:= f^{-1} (\s)$).

 Deduce that the fundamental classes  $E_i$ of the
 curves given by the irreducible components  in the fibers are linearly independent.

 Use Poincar\'e duality to deduce that the intersection form (cup product)
 matrix $||E_i \cdot E_j||$
 on these classes is nondegenerate.

 (Grauert proved a general theorem,
 valid in the analytic context and for an analytic  germ $(Y,o)$ that even shows that
 this form is negative definite.)

 Show that the contribution $\IH^*(Y)$ to $H^*(X)$ can be viewed
 as the space orthogonal, with respect to the cup product, to the span of the $E_i$'s.

 Deduce that $\IH^*(Y)$ sits inside $H^*(X,\rat)$ compatibly with the
 Hodge decomposition of $H^*(X, \comp)$, i.e., $\IH^j(Y,\rat)$ inherits a pure Hodge structure
 of weight $j$.
\index{Hodge structure}

 \end{enumerate}

\end{document}